# A SHORT NOTE FOR THE ROBUSTNESS PROPERTIES OF HYBRID DEAD-BEAT OBSERVERS

**Iasson Karafyllis*** and **Zhong-Ping Jiang****

*Dept. of Environmental Eng., Technical University of Crete,
73100, Chania, Greece, email: ikarafyl@enveng.tuc.gr

**Dept. of Electrical and Computer Eng., Polytechnic Institute of New York University,
Six Metrotech Center, Brooklyn, NY 11201, U.S.A., email: zjiang@control.poly.edu

**Abstract**

A discussion of the robustness properties of the proposed observer with respect to measurement errors is provided for the recently proposed full-order and reduced-order, hybrid, dead-beat observer for a class of nonlinear systems, linear in the unmeasured states.

**Keywords:** observer design, nonlinear systems, hybrid observers, robustness.

## 1. Introduction

In the recent paper [3] we studied the possibility of designing hybrid dead-beat observers for nonlinear systems of the form:

$$\dot{x}(t) = A(y(t), u(t))x(t) + b(y(t), u(t))$$
$$\dot{y}_i(t) = f_i(y(t), u(t)) + \sum_{j=1}^{n} c_{i,j}(y(t))x_j(t) \quad , \quad i = 1,...,k \tag{1.1}$$
$$x(t) = (x_1(t),...,x_n(t))' \in \Re^n , y(t) = (y_1(t),...,y_k(t))' \in \Re^k$$
$$(x(t), y(t)) \in O \subseteq \Re^n \times \Re^k , u(t) \in U \subseteq \Re^m$$

where $O \subseteq \Re^{n+k}$ is an open set, $U \subseteq \Re^m$ is a non-empty closed set, $A(y,u) = \{a_{i,j}(y,u), i,j = 1,...,n\}$ and all mappings $a_{i,j} : \Omega \times U \to \Re$ ($i, j = 1,...,n$), $b : \Omega \times U \to \Re^n$, $c_{i,j} : \Omega \to \Re$ ($i = 1,...,k$, $j = 1,...,n$) and $f_i : \Omega \times U \to \Re$ ($i = 1,...,k$) are locally Lipschitz, where $\Omega = \{y \in \Re^k : \exists x \text{ such that } (x, y) \in O\}$. It is assumed that the component of the state vector $y$, also known as the output, is available and that the remaining state component $x$ is unmeasured and is to be estimated.

The results of [3] were exploited in the recent work [4], where the design of hybrid dead-beat observers for chemostat models was studied.

The application dealing with the estimation of the frequency of a sinusoidal signal in [3] showed that the proposed hybrid dead-beat observer is robust with respect to high frequency noise. The results showed that the sensitivity to measurement noise decreases as the time horizon of the minimized $L^2$ norm increases, i.e., as the length of the history of the output which is utilized for the state estimation increases. This feature is expected and it is common to optimization-based observers. Motivated by these results, this short note is devoted to the study of the robustness properties of the proposed observer with respect to measurement errors. Proposition 3.1 implies that the difference of the state of system (1.1) and the observer state satisfies the Bounded-Input-Bounded-Output (BIBO) property (statement (a) of Proposition 3.1) and the Converging-Input-Converging-Output (CICO) property (statement (c) of Proposition 3.1) with the measurement error as input, under certain hypotheses. The result is important, because the topic of the robustness properties of observers for nonlinear systems is rarely studied by both numerical and theoretical tools.

**Notations** Throughout this note we adopt the following notations:

* Let $I \subseteq \mathfrak{R}_+ := [0,+\infty)$ be an interval. By $L^\infty(I;U)$ ($L^\infty_{loc}(I;U)$) we denote the space of measurable and (locally) essentially bounded functions $u(\cdot)$ defined on $I$ and taking values in $U \subseteq \mathfrak{R}^m$.

* By $C^0(A;\Omega)$, we denote the class of continuous functions on $A$, which take values in $\Omega$.

* For a vector $x \in \mathfrak{R}^n$ we denote by $x'$ its transpose and by $|x|$ its Euclidean norm. The determinant of a square matrix $A \in \mathfrak{R}^{n \times n}$ is denoted by $\det(A)$. $A' \in \mathfrak{R}^{n \times m}$ denotes the transpose of the matrix $A \in \mathfrak{R}^{m \times n}$.

## 2. Review of Hybrid Dead-Beat Observer Design

Consider an autonomous system described by ordinary differential equations of the form:

$$\dot{x}(t) = F(x(t), u(t))$$
$$x(t) \in D \subseteq \mathfrak{R}^n , u(t) \in U \subseteq \mathfrak{R}^m \tag{2.1}$$

where $D \subseteq \mathfrak{R}^n$ is an open set, $U \subseteq \mathfrak{R}^m$ is a non-empty closed set and the mapping $F: D \times U \to \mathfrak{R}^n$ is locally Lipschitz. The output of system (2.1) is given by

$$y(t) = h(x(t)) \tag{2.2}$$

where the mapping $h: D \to \mathfrak{R}^k$ is continuous. For system (2.1) we adopt the following notion of observability. We assume that for every $x_0 \in D$ and $u \in L^\infty_{loc}(\mathfrak{R}_+;U)$ there exists a unique solution $[0,+\infty) \ni t \to x(t) = x(t,x_0;u) \in D$ satisfying (2.1) for almost every $t \geq 0$ with $x(0) = x(0,x_0;u) = x_0$.

**Definition 2.1 and Definition 2.7 in [3]:** *Consider system (2.1) with output (2.2). We say that the input $u \in L^\infty([0,r];U)$ strongly distinguishes the state $x_0 \in D$ in time $r > 0$, if the following condition holds*

$$\max_{t \in [0,r]} |h(x(t,x_0;u)) - h(x(t,\xi;u))| > 0 \text{, for all } \xi \in D \text{ with } x_0 \neq \xi \tag{2.3}$$

*We say that (2.1) is strongly observable in time $r > 0$ if every input $u \in L^\infty([0,r];U)$ strongly distinguishes every state $x_0 \in D$ in time $r > 0$.*

For system (1.1) we assume that for every $(x_0, y_0) \in O$ and $u \in L^\infty_{loc}(\mathfrak{R}_+;U)$ there exists a unique mapping $[0,+\infty) \ni t \to (x(t), y(t)) = (x(t,x_0,y_0;u), y(t,x_0,y_0;u)) \in O$ satisfying (1.1) for almost every $t \geq 0$ with $(x(0), y(0)) = (x(0,x_0,y_0;u), y(0,x_0,y_0;u)) = (x_0, y_0)$.

We denote by $\Phi(t,x_0,y_0;u)$ the transition matrix of the linear time-varying system $\dot{x}(t) = A(y(t),u(t))x(t)$ when $u \in L^\infty_{loc}(\mathfrak{R}_+;U)$ and $y(t) = y(t,x_0,y_0;u)$ for $t \geq 0$ are considered as the inputs. Then the following fact holds for the solutions of system (1.1). It follows directly from integration of the differential equations (1.1).

Fact I: For every $(x_0, y_0) \in O$ and $u \in L^\infty([0,r];U)$ the following equations hold for all $t \geq 0$:

$$x(t,x_0,y_0;u) = \Phi(t,x_0,y_0;u)x_0 + \theta(t,x_0,y_0;u) \tag{2.4}$$

$$p(t,x_0,y_0;u) = q'(t,x_0,y_0;u) x_0 \tag{2.5}$$

where

$$q(t,x_0,y_0;u) := \int_0^t \Phi'(s,x_0,y_0;u)C(s,x_0,y_0;u)ds \tag{2.6}$$

$$\theta(t, x_0, y_0; u) := \int_0^t \Phi(t, x_0, y_0; u)\Phi^{-1}(\tau, x_0, y_0; u)b(y(\tau, x_0, y_0; u), u(\tau))d\tau \tag{2.7}$$

$$C'(t, x_0, y_0; u) := \begin{bmatrix} c_{1,1}(y(t, x_0, y_0; u)) & \cdots & c_{1,n}(y(t, x_0, y_0; u)) \\ \vdots & & \vdots \\ c_{k,1}(y(t, x_0, y_0; u)) & \cdots & c_{k,n}(y(t, x_0, y_0; u)) \end{bmatrix} \in \Re^{k \times n} \tag{2.8}$$

$$p(t, x_0, y_0; u) := y(t, x_0, y_0; u) - y_0 - \begin{bmatrix} \int_0^t f_1(y(s, x_0, y_0; u), u(s))ds \\ \vdots \\ \int_0^t f_k(y(s, x_0, y_0; u), u(s))ds \end{bmatrix} - \int_0^t C'(s, x_0, y_0; u)\theta(s, x_0, y_0; u)ds \tag{2.9}$$

It is important to note at this point that all expressions involved in (2.4)-(2.9) can be evaluated by means of the output trajectory $y(\tau) = y(\tau, x_0, y_0; u)$ for $\tau \in [0, t]$ and the input $u(\tau)$ for $\tau \in [0, t]$. For example, the transition matrix $\Phi(t, x_0, y_0; u)$ can be evaluated by solving the linear matrix differential equation $\frac{d}{d\tau}\Phi(\tau) = A(y(\tau), u(\tau))\Phi(\tau)$ for $\tau \in [0, t]$ with initial condition $\Phi(0) = I$, where $I$ denotes the identity matrix. Similarly, $C(\tau) := C(\tau, x_0, y_0; u)$ is simply $C'(\tau) := \{c_{i,j}(y(\tau)), i = 1, ..., k, j = 1, ..., n\}$ for $\tau \in [0, t]$ and $\theta(t) = \theta(t, x_0, y_0; u)$ can be computed by solving the linear system of differential equations $\frac{d}{d\tau}\theta(\tau) = A(y(\tau), u(\tau))\theta(\tau) + b(y(\tau), u(\tau))$ for $\tau \in [0, t]$ with initial condition $\theta(0) = 0 \in \Re^n$. Finally, the differential equations $\frac{d}{d\tau}q(\tau) = \Phi'(\tau)C(\tau)$ and $\frac{d}{d\tau}\xi(\tau) = (f_1(y(\tau), u(\tau)), ..., f_k(y(\tau), u(\tau)))' + C'(\tau)\theta(\tau)$, for $\tau \in [0, t]$ can be utilized to provide the quantities $q(t) = q(t, x_0, y_0; u)$ and $p(t, x_0, y_0; u) = y(t) - y(0) - \xi(t)$.

The following proposition provides characterizations of the class of inputs $u \in L^\infty([0, r]; U)$ which strongly distinguish the state $(x_0, y_0) \in O$ in time $r > 0$ for system (1.1). The basic idea is the conversion of the observability property to the minimization of an appropriate $L^2$ norm.

**Proposition 2.3 in [3]:** *Consider system (1.1). The following statements are equivalent:*

**(a)** *The input $u \in L^\infty([0, r]; U)$ strongly distinguishes the state $(x_0, y_0) \in O$ in time $r > 0$.*

**(b)** *The problem*

$$\min_{\xi \in B(y_0)} \int_0^r |p(t, x_0, y_0; u) - q'(t, x_0, y_0; u)\xi|^2 dt \tag{2.10}$$

*where $B(y_0) := \{\xi \in \Re^n : (\xi, y_0) \in O\}$, admits the unique solution $\xi = x_0$.*

**(c)** *The symmetric matrix*

$$Q(r, x_0, y_0; u) := \int_0^r q(t, x_0, y_0; u)q'(t, x_0, y_0; u)dt \tag{2.11}$$

*is positive definite. Moreover, it holds that*

$$x_0 = Q^{-1}(r, x_0, y_0; u) \int_0^r q(t, x_0, y_0; u) p(t, x_0, y_0; u) dt \qquad (2.12)$$

**(d)** *The following implication holds:*

$$q'(t, x_0, y_0; u) \xi = 0 \quad , \forall t \in [0, r] \Rightarrow \xi = 0 \in \Re^n \qquad (2.13)$$

The following corollary presents sufficient conditions for an input $u \in L^\infty([0,r]; U)$ to strongly distinguish $(x_0, y_0) \in O$ in time $r > 0$ for systems of the form (1.1) with scalar output.

**Corollary 2.5 in [3]:** *Consider system (1.1) with $k = 1$ and let $(x_0, y_0) \in O$, $u \in L^\infty([0,r]; U)$ for which there exist $t, t_1, \ldots, t_{n-1} \in [0, r]$ such that*

$$\det\left(\begin{bmatrix} C'(t, x_0, y_0; u)\Phi(t, x_0, y_0; u) \\ C'(t_1, x_0, y_0; u)\Phi(t_1, x_0, y_0; u) \\ \vdots \\ C'(t_{n-1}, x_0, y_0; u)\Phi(t_{n-1}, x_0, y_0; u) \end{bmatrix}\right) \neq 0 \qquad (2.15)$$

*Then the input $u \in L^\infty([0,r]; U)$ strongly distinguishes the state $(x_0, y_0) \in O$ in time $r > 0$. Moreover, the symmetric matrix $Q(r, x_0, y_0; u)$ defined by (2.11) is positive definite and (2.12) holds.*

Under the following hypothesis for system (1.1):

**(H1)** *System (1.1) is strongly observable in time $r > 0$.*

proposition 2.3 in [3] shows that we can define the operator:

$$P : L^\infty([0,r]; \Omega) \times L^\infty([0,r]; U) \to \Re^n$$

where $\Omega = \{y \in \Re^k : \exists x \text{ such that } (x, y) \in O\}$. For each $(y, u) \in L^\infty([0,r]; \Omega) \times L^\infty([0,r]; U)$, $P(y, u)$ is defined by

$$P(y, u) = \Phi(r, y; u) Q^{-1} \int_0^r q(\tau) p(\tau) d\tau + \theta(r) \qquad (2.16)$$

where $\Phi(t, y; u)$ is the transition matrix of the linear system $\dot{z}(t) = A(y(t), u(t)) z(t)$, $Q = \int_0^r q(\tau) q'(\tau) d\tau$, $q(\tau) = \int_0^\tau \Phi'(s, y; u) C(s) ds$, $C'(\tau) := \{c_{i,j}(y(\tau)), i = 1, \ldots, k, j = 1, \ldots, n\}$, $p(\tau) = y(\tau) - y(0) - \int_0^\tau f(y(s), u(s)) ds - \int_0^\tau C'(s) \theta(s) ds$, $f(y, u) := (f_1(y, u), \ldots, f_k(y, u))'$, $\theta(\tau) := \int_0^\tau \Phi(\tau, y; u) \Phi^{-1}(s, y; u) b(y(s), u(s)) ds$ for all $\tau \in [0, r]$. Proposition 2.3 in [3] guarantees that, if hypothesis (H1) holds for system (1.1), then for every $(x_0, y_0) \in O$ and $u \in L^\infty_{loc}(\Re_+; U)$ the following equality holds:

$$x(t, x_0, y_0; u) = P(\delta_{t-r} y, \delta_{t-r} u), \text{ for all } t \geq r \qquad (2.17)$$

where $(\delta_{t-r} y)(s) = y(t - r + s, x_0, y_0; u)$, $(\delta_{t-r} u)(s) = u(t - r + s)$ for $s \in [0, r]$.

Therefore, if hypothesis (H1) holds for system (1.1), then we are in a position to provide a hybrid, dead-beat observer for system (1.1). Given $t_0 \geq 0$, $(z_0, w_0) \in O$, we calculate $(z(t), w(t)) \in O$ by the following algorithm:

Step $i$: Calculation of $z(t)$ for $t \in [t_0 + ir, t_0 + (i+1)r]$

1) Calculate $z(t)$ for $t \in [t_0 + ir, t_0 + (i+1)r)$ as the solution of $\dot{z}(t) = A(w(t), u(t))z(t) + b(w(t), u(t))$, $\dot{w}_i(t) = f_i(w(t), u(t)) + \sum_{j=1}^{n} c_{i,j}(w(t))z_j(t)$ ($i = 1, ..., k$), with $w(t) = (w_1(t), ..., w_k(t))' \in \Re^k$.

2) Set $z(t_0 + (i+1)r) = P(\delta_{t_0 + ir} y, \delta_{t_0 + ir} u)$ and $w(t_0 + (i+1)r) = y(t_0 + (i+1)r)$, where $P: L^\infty([0, r]; \Omega) \times L^\infty([0, r]; U) \to \Re^n$ is the operator defined by (2.16).

For $i = 0$ we take $(z(t_0), w(t_0)) = (z_0, w_0)$ (initial condition).

The proposed observer can be represented by the following system of equations:

$$\begin{aligned}
&\dot{z}(t) = A(w(t), u(t))z(t) + b(w(t), u(t)), \; t \in [\tau_i, \tau_{i+1}) \\
&\dot{w}_i(t) = f_i(w(t), u(t)) + \sum_{j=1}^{n} c_{i,j}(w(t))x_j(t), \quad i = 1, ..., k, \; t \in [\tau_i, \tau_{i+1}) \\
&z(\tau_{i+1}) = P(\delta_{\tau_i} y, \delta_{\tau_i} u) \\
&w(\tau_{i+1}) = y(\tau_{i+1}) \\
&\tau_{i+1} = \tau_i + r \\
&z(t) = (z_1(t), ..., z_n(t))' \in \Re^n, \; w(t) = (w_1(t), ..., w_k(t))' \in \Re^k \\
&(z(t), w(t)) \in O \subseteq \Re^n \times \Re^k
\end{aligned} \quad (2.18)$$

Thus, from all the above results, we obtain the following corollary.

**Corollary 2.9 in [3]:** *Consider system (1.1) and assume that hypothesis (H1) holds. Consider the unique solution $(x(t), y(t), z(t), w(t)) \in O \times O$ of (1.1), (2.18) with arbitrary initial condition $(x_0, y_0, z_0, w_0) \in O \times O$ corresponding to arbitrary input $u \in L^\infty_{loc}(\Re_+; U)$. Then the solution $(x(t), y(t), z(t), w(t)) \in O \times O$ of (1.1), (2.18) satisfies:*

$$z(t) = x(t) \text{ and } w(t) = y(t), \text{ for all } t \geq r \quad (2.19)$$

Next assume that the following hypothesis holds in addition to hypothesis (H1).

**(H2)** *There exist open sets $D \subseteq \Re^n$ and $\Omega \subseteq \Re^k$ such that $O = D \times \Omega$. Moreover, for every $\xi \in D$ and for every $(y, u) \in L^\infty([0, r]; \Omega) \times L^\infty([0, r]; U)$, the solution $z(t)$ of $\dot{z}(t) = A(y(t), u(t))z(t) + b(y(t), u(t))$ with $z(0) = \xi$ satisfies $z(t) \in D$ for all $t \in [0, r]$.*

If hypothesis (H2) holds, then we can design a reduced-order, hybrid, dead-beat observer for system (1.1) of the form:

$$\begin{aligned}
&\dot{z}(t) = A(y(t), u(t))z(t) + b(y(t), u(t)), \; t \in [\tau_i, \tau_{i+1}) \\
&z(\tau_{i+1}) = P(\delta_{\tau_i} y, \delta_{\tau_i} u) \\
&\tau_{i+1} = \tau_i + r \\
&z(t) = (z_1(t), ..., z_n(t))' \in D \subseteq \Re^n
\end{aligned} \quad (2.20)$$

where $P: L^\infty([0, r]; \Omega) \times L^\infty([0, r]; U) \to \Re^n$ is the operator defined by (2.16).

**Corollary 2.11 in [3]:** *Consider system (1.1) and assume that hypotheses (H1), (H2) hold. Consider the unique solution $(x(t), y(t), z(t)) \in D \times \Omega \times D$ of (1.1), (2.20) with arbitrary initial condition $(x_0, y_0, z_0) \in D \times \Omega \times D$ corresponding to arbitrary input $u \in L^\infty_{loc}(\Re_+; U)$. Then the solution $(x(t), y(t), z(t)) \in D \times \Omega \times D$ of (1.1), (2.20) satisfies:*

$$z(t) = x(t), \text{ for all } t \geq r \qquad (2.21)$$

## 3. Robustness Issues

In this section, we discuss the robustness properties of the proposed observer. We focus on the case where both hypotheses (H1), (H2) are satisfied with $D = \Re^n$, $\Omega = \Re^k$ (although the same discussion can be applied to systems satisfying only hypothesis (H1)). It should be noted that a systematic study of the robustness properties of observers for nonlinear systems is rare and the topic is completely "untouched".

Most specifically, we study the robustness properties of the proposed observer subject to measurement errors. By measurement error we mean a measurable and locally essentially bounded input $e: \Re_+ \to \Re^k$ which corrupts the output values that are fed to the observer, i.e., the observer is described by the equations:

$$\begin{aligned}
\dot{z}(t) &= A(\tilde{y}(t), u(t))z(t) + b(\tilde{y}(t), u(t)), \ t \in [\tau_i, \tau_{i+1}) \\
z(\tau_{i+1}) &= P(\delta_{\tau_i}\tilde{y}, \delta_{\tau_i}u) \\
\tau_{i+1} &= \tau_i + r \\
z(t) &= (z_1(t), \ldots, z_n(t))' \in \Re^n
\end{aligned} \qquad (3.1)$$

where $P: L^\infty([0,r]; \Re^k) \times L^\infty([0,r]; U) \to \Re^n$ is the operator defined by (2.16) and

$$\tilde{y}(t) = y(t) + e(t), \ \forall t \geq 0 \qquad (3.2)$$

Taking into account the formula $z(\tau_{i+1}) = P(\delta_{\tau_i}\tilde{y}, \delta_{\tau_i}u) = P(\delta_{\tau_i}y + \delta_{\tau_i}e, \delta_{\tau_i}u)$ for all $i \geq 0$ and Corollary 2.11 in [3], which guarantees that $x(\tau_{i+1}) = P(\delta_{\tau_i}y, \delta_{\tau_i}u)$ for all $i \geq 0$, we can conclude that the observer error induced by the measurement error at $t = \tau_{i+1}$ ($i \geq 0$) will satisfy:

$$|z(\tau_{i+1}) - x(\tau_{i+1})| = |P(\delta_{\tau_i}y + \delta_{\tau_i}e, \delta_{\tau_i}u) - P(\delta_{\tau_i}y, \delta_{\tau_i}u)|, \text{ for all } i \geq 0 \qquad (3.3)$$

At this point, a pair of hypotheses is introduced in order to analyze further the time evolution of the observer error.

**(R1)** *For every $(x_0, y_0) \in \Re^n \times \Re^k$ and $u \in L^\infty(\Re_+; U)$, the unique solution $(x(t), y(t)) \in \Re^n \times \Re^k$ of (1.1), with initial condition $(x_0, y_0) \in \Re^n \times \Re^k$ corresponding to input $u \in L^\infty(\Re_+; U)$ satisfies $\sup_{t \geq 0}|y(t)| + \sup_{t \geq 0}|x(t)| < +\infty$.*

**(R2)** *The operator $L^\infty([0,r]; \Re^k) \times L^\infty([0,r]; U) \ni (y, u) \to P(y, u) \in \Re^n$ is completely continuous with respect to $y \in L^\infty([0,r]; \Re^k)$, i.e., for every pair of bounded sets $S \subset L^\infty([0,r]; \Re^k)$, $V \subseteq L^\infty([0,r]; U)$ the image set $P(S \times V) \subset \Re^n$ is bounded and for every $\varepsilon > 0$ there exists $\delta > 0$ such that $|P(y, u) - P(\hat{y}, u)| < \varepsilon$ for every $y, \hat{y} \in S$, $u \in V$ with $\sup_{t \in [0,r]}|y(t) - \hat{y}(t)| < \delta$.*

Hypothesis (R1) imposes restrictions on the dynamic behavior of system (1.1). On the other hand, hypothesis (R2) is a continuity hypothesis which can be guaranteed easily for certain cases. A case where hypothesis (R2) holds is the case where for every pair of bounded sets $S \subset L^\infty([0,r]; \Re^k)$, $V \subseteq L^\infty([0,r]; U)$ there exists $a > 0$ such that $\det(Q) \geq a$ for all $y \in S$, $u \in V$, where $Q = \int_0^r q(\tau)q'(\tau)d\tau$, $q(\tau) = \int_0^\tau \Phi'(s,y;u)C(s)ds$, $C'(\tau) := \{c_{i,j}(y(\tau)), i=1,...,k, j=1,...,n\}$ and $\Phi(t,y;u)$ is the transition matrix of the linear system $\dot{z}(t) = A(y(t),u(t))z(t)$. This reminds the case of uniform observability of linear time-varying systems (see [2]).

Another thing that should be noted here is that the estimation of the state of system (1.1) under hypothesis (R1) cannot be performed in general by means of a high-gain observer (see [1]). Indeed, although the subsystem $\dot{x}(t) = A(y(t),u(t))x(t) + b(y(t),u(t))$ is globally Lipschitz when (R1) holds and a bounded input $u \in L^\infty(\Re_+;U)$ is applied, we are not aware of the Lipschitz constant of the system (since we do not assume knowledge of $\sup_{t\geq 0}|y(t)| < +\infty$), which in general will depend on the initial conditions $(x_0, y_0) \in \Re^n \times \Re^k$ and the applied input $u \in L^\infty(\Re_+;U)$.

Using hypotheses (R1), (R2), we are in a position to show the following robustness result.

**Proposition 3.1:** *Consider system (1.1) for which hypotheses (H1), (H2) hold with $D = \Re^n$, $\Omega = \Re^k$. Moreover, assume that hypotheses (R1), (R2) hold as well. Then*

a) *for every $(x_0, y_0, z_0) \in \Re^n \times \Re^k \times \Re^n$, $(u,e) \in L^\infty(\Re_+;U) \times L^\infty(\Re_+;\Re^k)$ the solution $(x(t), y(t), z(t)) \in \Re^n \times \Re^k \times \Re^n$ of (1.1), (3.1), (3.2) with initial condition $(x_0, y_0, z_0) \in \Re^n \times \Re^k \times \Re^n$ corresponding to inputs $(u,e) \in L^\infty(\Re_+;U) \times L^\infty(\Re_+;\Re^k)$ satisfies $\sup_{t\geq 0}|z(t) - x(t)| < +\infty$,*

b) *for every $(x_0, y_0, z_0) \in \Re^n \times \Re^k \times \Re^n$, $u \in L^\infty(\Re_+;U)$ and $\varepsilon > 0$ there exists $\delta > 0$ such that for every $e \in L^\infty(\Re_+;\Re^k)$ with $\sup_{t\geq 0}|e(t)| < \delta$ the solution $(x(t), y(t), z(t)) \in \Re^n \times \Re^k \times \Re^n$ of (1.1), (3.1), (3.2) with initial condition $(x_0, y_0, z_0) \in \Re^n \times \Re^k \times \Re^n$ corresponding to inputs $(u,e) \in L^\infty(\Re_+;U) \times L^\infty(\Re_+;\Re^k)$ satisfies $\sup_{t\geq r}|z(t) - x(t)| < \varepsilon$,*

c) *for every $(x_0, y_0, z_0) \in \Re^n \times \Re^k \times \Re^n$, $(u,e) \in L^\infty(\Re_+;U) \times L^\infty(\Re_+;\Re^k)$ with $\lim_{t\to+\infty}|e(t)| = 0$, the solution $(x(t), y(t), z(t)) \in \Re^n \times \Re^k \times \Re^n$ of (1.1), (3.1), (3.2) with initial condition $(x_0, y_0, z_0) \in \Re^n \times \Re^k \times \Re^n$ corresponding to inputs $(u,e) \in L^\infty(\Re_+;U) \times L^\infty(\Re_+;\Re^k)$ satisfies $\lim_{t\to+\infty}|z(t) - x(t)| = 0$.*

**Proof:** Let $(x_0, y_0, z_0) \in \Re^n \times \Re^k \times \Re^n$, $(u,e) \in L^\infty(\Re_+;U) \times L^\infty(\Re_+;\Re^k)$ be given and consider the solution $(x(t), y(t), z(t)) \in \Re^n \times \Re^k \times \Re^n$ of (1.1), (3.1), (3.2) with initial condition $(x_0, y_0, z_0) \in \Re^n \times \Re^k \times \Re^n$ corresponding to inputs $(u,e) \in L^\infty(\Re_+;U) \times L^\infty(\Re_+;\Re^k)$. By virtue of hypothesis (R1) we have $\|y\| = \sup_{t\geq 0}|y(t)| < +\infty$ and $\|x\| = \sup_{t\geq 0}|x(t)| < +\infty$. The proof is based on the following fact, which exploits the fact that $A(y,u) = \{a_{i,j}(y,u), i,j=1,...,n\}$ and all mappings $a_{i,j}: \Re^k \times U \to \Re$ ($i,j=1,...,n$), $b: \Re^k \times U \to \Re^n$, are locally Lipschitz. Its proof is standard and is omitted.

**Fact:** There exist non-decreasing functions $\kappa: \Re_+ \to \Re_+$ and $R: \Re_+ \to \Re_+$ such that

$$|x(t) - z(t)| \leq \left[|x(\tau_i) - z(\tau_i)| + r\|e\|R(s)\right]\exp(r\kappa(s)), \quad \forall t \in (\tau_i, \tau_{i+1}) \tag{3.4}$$

where $s := \|y\| + \|e\| + \|u\| + \|x\|$, $\|y\| = \sup_{t\geq 0}|y(t)| < +\infty$, $\|x\| = \sup_{t\geq 0}|x(t)| < +\infty$, $\|u\| = \sup_{t\geq 0}|u(t)| < +\infty$ and $\|e\| = \sup_{t\geq 0}|e(t)| < +\infty$.

By virtue of hypothesis (R2), which implies that the image set $P(S \times V) \subset \Re^n$ is bounded, where $S \subset L^\infty([0,r];\Re^k)$ is the bounded set of measurable and essentially bounded functions $z:[0,r] \to \Re^k$ with $\sup_{t \geq 0}|z(t)| \leq \|y\| + \|e\|$ and $V \subseteq L^\infty([0,r];U)$ is the bounded set of measurable and essentially bounded functions $v:[0,r] \to U$ with $\sup_{t \geq 0}|v(t)| \leq \|u\|$, there exists $K \geq 0$ such that $|P(\delta_{\tau_i} y + \delta_{\tau_i} e, \delta_{\tau_i} u)| \leq K$ and $|P(\delta_{\tau_i} y, \delta_{\tau_i} u)| \leq K$, for all $i \geq 0$. Consequently, we obtain from (3.3):

$$|z(\tau_{i+1}) - x(\tau_{i+1})| \leq 2K, \text{ for all } i \geq 0 \qquad (3.5)$$

Combining (3.4) and (3.5) we can conclude that $\sup_{t \geq 0}|z(t) - x(t)| < +\infty$.

Inequality (3.4) implies that for every $\varepsilon > 0$ there exists $\delta_1 \in (0, \varepsilon)$ such that $|x(t) - z(t)| < \varepsilon$ for all $t \in (\tau_i, \tau_{i+1})$, provided that $\|e\| < \delta_1$ and $|z(\tau_i) - x(\tau_i)| < \delta_1$. Moreover, hypothesis (R2) and (3.3) implies the existence of $\delta_2 \in (0, \delta_1)$ such that $|z(\tau_{i+1}) - x(\tau_{i+1})| < \delta_1$, for all $i \geq 0$, provided that $\|e\| < \delta_2$. Combining the two previous inequalities, we conclude that $|x(t) - z(t)| < \varepsilon$ for all $t \geq r$, provided that $\|e\| < \min\{\delta_1, \delta_2\}$. On the other hand, if $\lim_{t \to +\infty}|e(t)| = 0$ then there exists $i \geq 0$ sufficiently large such that $\sup_{t \geq \tau_i}|e(t)| < \min\{\delta_1, \delta_2\}$. The previous inequalities imply that $|x(t) - z(t)| < \varepsilon$ for all $t \geq \tau_i + r$ and since $\varepsilon > 0$ is arbitrary we conclude that $\lim_{t \to +\infty}|z(t) - x(t)| = 0$.

The proof is complete. ◁

**Remark 3.2:** It should be noted that if $A(y,u) = A(u)$, then there is no need to assume that $\sup_{t \geq 0}|x(t)| < +\infty$. Therefore, in this case the conclusions of Proposition 3.1 hold if hypothesis (R1) is replaced by the following weaker hypothesis:

**(R3)** *For every $(x_0, y_0) \in \Re^n \times \Re^k$ and $u \in L^\infty(\Re_+; U)$, the unique solution $(x(t), y(t)) \in \Re^n \times \Re^k$ of (1.1), with initial condition $(x_0, y_0) \in \Re^n \times \Re^k$ corresponding to input $u \in L^\infty(\Re_+; U)$ satisfies $\sup_{t \geq 0}|y(t)| < +\infty$.*

**Example 3.3:** Consider the system

$$\begin{aligned} \dot{x}(t) &= a(y(t), u(t))x(t) \\ \dot{y}(t) &= f(y(t), u(t)) + c(y(t))x(t) \\ x(t) &\in D, \ y(t) \in \Omega, \ u(t) \in U \end{aligned} \qquad (3.6)$$

with $D = \Re$, $\Omega = \Re$, which was studied in Example 2.12 in [3]. Here, we assume that hypothesis (R1) holds for system (3.6) and that $c(y) > 0$ for all $y \in \Re$. Therefore, hypothesis (H3) in [3] holds automatically and the reader can verify that hypothesis (R2) holds for the mapping $P(y,u)$ defined by

$$P(y,u) = \exp\left(\int_0^r a(y(s), u(s))ds\right) \frac{\int_0^r \left(y(\tau) - y(0) - \int_0^\tau f(y(s), u(s))ds\right)\left(\int_0^\tau c(y(s))\exp\left(\int_0^s a(y(w), u(w))dw\right)ds\right)d\tau}{\int_0^r \left(\int_0^\tau c(y(s))\exp\left(\int_0^s a(y(w), u(w))dw\right)ds\right)^2 d\tau}$$

(3.7)

Therefore, Proposition 3.1 guarantees the BIBO (Bounded-Input-Bounded-Output) and CICO (Converging-Input-Converging-Output) properties for the output $Y(t) = z(t) - x(t)$ from the input $e \in L^\infty(\Re_+; \Re)$ for system (3.6) and the dead-beat hybrid observer

$$\dot{z}(t) = a(y(t)+e(t),u(t))z(t), \; t \in [\tau_i, \tau_{i+1})$$
$$z(\tau_{i+1}) = P(\delta_{\tau_i} y + \delta_{\tau_i} e, \delta_{\tau_i} u) \quad (3.8)$$
$$\tau_{i+1} = \tau_i + r$$

where $P(y,u)$ is defined by (3.7). ◁